\documentstyle[12pt]{article}

\title{The stability of magnetic vortices\footnote{Research on this paper
       was supported by NSERC under grant N7901}}
       
\author{S. Gustafson  \hspace{1cm} I.M. Sigal\\
        gustaf@math.toronto.edu \hspace{0.5cm}  
        sigal@math.toronto.edu \\
        Dept. of Mathematics, University of Toronto \\
        100 St. George St., Toronto, ON, Canada.  M5S 3G3 \\
        Fax: (416) 978-4107}

\date{Dec 1, 1998}

\newtheorem{lemma}{Lemma}
\newtheorem{thm}{Theorem}
\newtheorem{rem}{Remark}
\newtheorem{prop}{Proposition}

\begin{document}

\maketitle

\abstract{We study the linearized stability of $n$-vortex ($n \in {\bf Z}$)
solutions of the magnetic Ginzburg-Landau (or Abelian Higgs) equations. 
We prove that the fundamental vortices ($n = \pm 1$) are stable for all
values of the coupling constant, $\lambda$, and
we prove that the higher-degree vortices ($|n| \geq 2$)
are stable for $\lambda < 1$, and unstable for $\lambda > 1$.
This resolves a long-standing conjecture (see, eg, \cite{jt}).}

\tableofcontents

\section{Introduction}

In this paper, we determine the stability of magnetic
(or Abelian Higgs) vortices.  
These are certain critical points of the energy functional 
\begin{equation}
\label{eq:ac}
E(\psi,A) = \frac{1}{2} \int_{{\bf R}^2} \left\{ |\nabla_A \psi|^2
                  + (\nabla \times A)^2 + \frac{\lambda}{4}
                  (|\psi|^2-1)^2 \right\}
\end{equation}
for the fields
\[  A : {\bf R}^2 \rightarrow {\bf R}^2 \;\;\;\;\;
\mbox{  and  }  \;\;\;\;\;
\psi : {\bf R}^2 \rightarrow {\bf C}. \] 
Here 
$\nabla_A = \nabla - iA$
is the covariant gradient, and $\lambda > 0$ is a coupling
constant.  For a vector, $A$, $\nabla \times A$ is the 
scalar $\partial_1 A_2 - \partial_2 A_1$, and for a scalar $\xi$,
$\nabla \times \xi$ is the vector $(-\partial_2 \xi, \partial_1 \xi)$.
Critical points of $E(\psi,A)$ satisfy the {\em Ginzburg-Landau}
(GL) equations
\begin{equation}
\label{eq:eq1}
  -\Delta_A\psi +
  \frac{\lambda}{2}(|\psi|^2-1)\psi = 0
\end{equation}
\begin{equation}
\label{eq:eq2}
  \nabla \times \nabla \times A
  - \Im(\bar{\psi} \nabla_A \psi) = 0
\end{equation}
where $\Delta_A = \nabla_A \cdot \nabla_A$.

Physically, the functional $E(\psi,A)$ gives the 
difference in free energy between the superconducting and normal states
near the transition temperature in the Ginzburg-Landau theory.
$A$ is the vector potential ($\nabla \times A$ is the induced
magnetic field), and $\psi$ is an {\em order parameter}.  The modulus
of $\psi$ is interpreted as describing the local density of 
superconducting Cooper pairs of electrons.

The functional $E(\psi,A)$ also gives the energy of a 
static configuration in the
Yang-Mills-Higgs classical gauge theory on ${\bf R}^2$, with abelian 
gauge group $U(1)$.  In this case $A$ is a connection
on the principal $U(1)$- bundle ${\bf R}^2 \times U(1)$, 
and $\psi$ is the {\em Higgs field} (see \cite{jt} for details).

A central feature of the functional $E(\psi,A)$ (and the GL equations)
is its infinite-dimensional symmetry group.  Specifically,
$E(\psi,A)$ is invariant under $U(1)$ {\em gauge transformations},
\begin{equation}
\label{eq:g1}
  \psi \mapsto e^{i\gamma}\psi
\end{equation}
\begin{equation}
\label{eq:g2}
  A \mapsto A + \nabla \gamma
\end{equation}
for any smooth $\gamma : {\bf R}^2 \rightarrow {\bf R}$.
In addition, $E(\psi,A)$ is invariant under coordinate
translations, and under the coordinate rotation transformation
\begin{equation}
\label{eq:rot}
  \psi(x) \mapsto \psi(g^{-1}x)  \;\;\;\;\;\;\;\;\;\; 
  A(x) \mapsto gA(g^{-1}x)
\end{equation}
for $g \in SO(2)$.

Finite energy field configurations satisfy
\begin{equation}
\label{eq:bc}
  |\psi| \rightarrow 1  \;\;\;\;\; \mbox{ as }
  \;\;\;\;\; |x| \rightarrow \infty
\end{equation}
which leads to the definition of the {\em topological degree},
$\mbox{deg}(\psi)$, of such a configuration:
\[
  \mbox{deg}(\psi) = \mbox{deg} \left(
  \left. \frac{\psi}{|\psi|} \right|_{|x| = R} : 
  {\bf S}^1 \rightarrow {\bf S}^1  \right)
\]
($R$ sufficiently large).  The degree is related to the phenomenon
of flux quantization.  Indeed, an application of Stokes' theorem shows that
a finite-energy configuration satisfies
\[
  \mbox{deg}(\psi) = \frac{1}{2\pi} \int_{{\bf R}^2} (\nabla \times A). 
\]

We study, in particular, ``radially-symmetric'' or ``equivariant''
fields of the form 
\begin{equation}
  \label{eq:psiA}
  \psi^{(n)}(x) = f_n(r)e^{in\theta}  \;\;\;\;\;\;\;\;\;\;
  A^{(n)}(x) = n\frac{a_n(r)}{r} \hat{x}^{\perp}
\end{equation}
where
$(r,\theta)$ are polar coordinates on ${\bf R}^2$,
$\hat{x}^{\perp} = \frac{1}{r} (-x_2,x_1)^t$,
$n$ is an integer, and  
\[
  f_n, a_n : [0,\infty) \rightarrow {\bf R}.
\]
It is easily checked that such configurations 
(if they satisfy~(\ref{eq:bc})) have degree $n$.
The existence of critical points of this form is well-known
(see section~\ref{subsec:vort}).  They are called {\em $n$-vortices}.

Our main results concern the stability of these
$n$-vortex solutions.  Let 
\[
  L^{(n)} = \mbox{ Hess } E (\psi^{(n)}, A^{(n)})
\]
be the linearized operator for GL around the $n$-vortex, acting on 
the space
\[
  X = L^2({\bf R}^2,{\bf C}) \oplus L^2({\bf R}^2,{\bf R}^2).
\]
The symmetry group of $E(\psi,A)$ gives rise to an infinite-dimensional
subspace of $\ker(L^{(n)}) \subset X$ (see section~\ref{subsec:zeromodes}),
which we denote here by $Z_{sym}$.  
We say the $n$-vortex is (linearly) {\em stable} if for some $c > 0$,
\[
  L^{(n)}|_{Z_{sym}^{\perp}} \geq c,
\]
and {\em unstable} if $L^{(n)}$ has a negative eigenvalue.  
The basic result of this paper is 
the following linearized stability statement:
\begin{thm}
\label{thm:main}
\begin{enumerate}
\item (Stability of fundamental vortices) \\
For all $\lambda > 0$, the $\pm1$-vortex is stable.
\item (Stability/instability of higher-degree vortices) \\
For $|n| \geq 2$, the $n$-vortex is
\[
  \left\{  \begin{array}{cc}
  \mbox{ stable }  & \mbox{ for } \lambda < 1  \\
  \mbox{ unstable }  & \mbox{ for } \lambda > 1.
  \end{array} \right.
\]
\end{enumerate}
\end{thm}

Theorem~\ref{thm:main} is the basic ingredient in a proof of the nonlinear
dynamical stability/instability of the $n$-vortex for certain 
dynamical versions of the GL equations.  
These include the GL gradient flow equations,
the Abelian Higgs (Lorentz-invariant) equations, and 
the Maxwell equations coupled to a nonlinear Schr\"odinger equation. 
These dynamical stability results are  established in a 
companion paper~(\cite{g2}).

The statement of theorem~\ref{thm:main} was conjectured in~\cite{jt}
on the basis of numerical observations (see~\cite{jr}).
Bogomolnyi (\cite{bog}) gave an argument for instability of vortices
for $\lambda > 1$, $|n| \geq 2$.  Our result
rigorously establishes this property.

The solutions of (\ref{eq:eq1}-\ref{eq:eq2}) are well-understood
in the case of {\em critical coupling}, $\lambda = 1$.
In this case, the {\em Bogomolnyi method} (\cite{bog}) gives a pair
of first-order equations whose solutions are global minimizers
of $E(\psi,A)$ among fields of fixed degree (and hence solutions of the
the GL equations).  Taubes (\cite{t1,t2})
has shown that all solutions of GL with $\lambda = 1$ are solutions
of these first-order equations, and that for a given degree $n$,
the gauge-inequivalent solutions 
form a $2|n|$-parameter family.
The $2|n|$ parameters describe the locations of the zeros of the scalar
field.  This is discussed in more detail in \cite{jt} 
(see also \cite{bgp}) and section~\ref{sec:crit}.  
We remark that for $\lambda = 1$, an $n$-vortex solution~(\ref{eq:psiA}) 
corresponds to the case when all $|n|$ zeros of the
scalar field lie at the origin.

The remainder of this paper is organized as follows.
In section~\ref{sec:vort} we describe in detail various properties
of the $n$-vortex.  In particular, we establish an important estimate
on the $n$-vortex profiles which differentiates between the cases
$\lambda < 1$ and $\lambda > 1$.
In section~\ref{sec:lin}, we introduce the linearized operator, fix the
gauge on the space of perturbations, and identify the zero-modes due to
symmetry-breaking.
Sections~\ref{sec:block} through \ref{sec:high} comprise a proof of 
theorem~\ref{thm:main}.  A block-decomposition for the linearized operator is
described in section~\ref{sec:block}.  This approach is similar to that
used to study the stability of non-magnetic vortices in~\cite{os1}
and~\cite{g1}.  In section~\ref{sec:fund},
we establish the positivity of certain blocks
(those corresponding to the radially-symmetric variational problem,
and those containing the translational zero-modes) for all $\lambda$,
which completes the stability proof for the $\pm 1$-vortices.  The basic
techniques are the characterization of symmetry-breaking in terms of
zero-modes of the Hessian (or linearized operator), and
a Perron-Frobenius type argument, 
based on a version of the maximum principle for systems
(proposition~\ref{prop:mainmp}), which shows that 
the translational zero-modes correspond to the bottom of the spectrum
of the linearized operator.  
A more careful analysis is needed for $|n| \geq 2$.  
This requires us to review some aspects of the critical case 
($\lambda = 1$) in section~\ref{sec:crit}.
The stability/instability proof for $|n| \geq 2$ is completed in
section~\ref{sec:high}.  We use an extension of Bogomolnyi's instability
argument, and another application of the Perron-Frobenius theory.

  


{\em Acknowledgment:  the first author would like to thank
the Courant institute for its hospitality during part of the
preparation of this paper, and especially J. Shatah for some
helpful discussions.  Part of this work is toward fulfillment of the
requirements of the first author's PhD at the University of Toronto.
The second author thanks Yu. N. Ovchinnikov
for many fruitful discussions.}

\section{The $n$-vortex}
\label{sec:vort}
 
In this section we discuss the existence, and properties,
of $n$-vortex solutions.

\subsection{Vortex solutions}
\label{subsec:vort}

The existence of solutions of (GL) of the form~(\ref{eq:psiA})
is well-known:

\begin{thm}[Vortex Existence; \cite{p,bc}]
\label{thm:exist}
  For every integer $n$, there is a solution 
  \begin{equation}
  \label{eq:nvortex}
    \psi^{(n)}(x) = f_n(r)e^{in\theta}  \;\;\;\;\;\;\;\;\;\;
    A^{(n)}(x) = n\frac{a_n(r)}{r} \hat{x}^{\perp}
  \end{equation}
  of the variational equations (\ref{eq:eq1})-(\ref{eq:eq2}).
  In particular, the radial functions ($f_n$, $a_n$)
  minimize the radial energy functional
  \begin{equation}
  \label{eq:erad}
    E^{(n)}_r(f,a) = \frac{1}{2} \int_0^{\infty}
    \left\{ (f')^2 + n^2\frac{(1-a)^2 f^2}{r^2}
    + n^2\frac{(a')^2}{r^2} + \frac{\lambda}{4}(f^2-1)^2 \right\} rdr
  \end{equation}
  (which is the full energy functional~(\ref{eq:ac}) restricted
  to fields of the form~(\ref{eq:psiA}))
  in the class
  \[
    \{ f,a : [0,\infty) \rightarrow {\bf R} \;\; | \;\;
    1-f \in H_1(rdr), \frac{a}{r} \in L^2_{loc}(rdr),
    \frac{a'}{r} \in L^2(rdr) \}.
  \]
  The functions $f_n$, $a_n$ are smooth, and have the 
  following properties (for $n \not= 0$):
  \begin{enumerate}
  \item $0 < f_n < 1$, $0 < a_n < 1$ on $(0,\infty)$
  \item $f_n', a_n' > 0$
  \item $f_n \sim cr^n$, $a_n \sim dr^2$,
        as $r \rightarrow 0$ ($c>0$ and $d>0$ are constants)
  \item $1-f_n$, $1-a_n \rightarrow 0$
        as $r \rightarrow \infty$, with an exponential rate of decay.
  \end{enumerate}
\end{thm}
We call ($\psi^{(n)}$, $A^{(n)}$) an {\em $n$-vortex} 
(centred at the origin).


It follows immediately that the functions $f_n$ and $a_n$ 
satisfy the ODEs
\begin{equation}
\label{eq:ode1}
  -\Delta_r f_n + \frac{n^2(1-a_n)^2}{r^2}f_n + 
  \frac{\lambda}{2}(f_n^2-1)f_n = 0 
\end{equation}
and
\begin{equation}
\label{eq:ode2}
  -a_n'' + \frac{a_n'}{r} - f_n^2(1-a_n) = 0.
\end{equation}
\begin{rem}
To our knowledge, it is not known if solutions of the form~(\ref{eq:psiA})
are unique.  In the appendix, we show that for $\lambda \geq 2n^2$,
any such solution minimizes $E_r^{(n)}$.
\end{rem}
\begin{rem}
The functions $f_n$ and $a_n$ also depend on $\lambda$, but
we suppress this dependence for ease of notation.  
When it will cause no confusion, we will also drop the subscript $n$.
\end{rem}
\begin{rem}
The discrete symmetry $\psi \mapsto \bar{\psi}$, $A \mapsto -A$ of (GL)
interchanges $(\psi^{(n)}, A^{(n)})$ and $(\psi^{(-n)}, A^{(-n)})$.
Thus, we can assume $n \geq 0$.
\end{rem}

\subsection{An estimate on the vortex profiles}


The following inequality, relating the exponentially decaying
quantities $f'$ and $1-a$, plays a crucial role in the
stability/instability proof.
\begin{prop}
We have
\begin{equation}
\label{eq:ineq}
  \left\{  \begin{array}{cc}
  f'(r) > \frac{n(1-a(r))}{r}f(r)  &  \;\;\; 
  \mbox{ for } \;\;\; \lambda < 1  \\
  f'(r) < \frac{n(1-a(r))}{r}f(r)  &  \;\;\; 
  \mbox{ for } \;\;\; \lambda > 1
  \end{array} \right.
\end{equation}
\end{prop}
{\em Proof:}
Define $e(r) \equiv f'(r) - \frac{n(1-a(r))}{r}f(r)$.
The properties listed in theorem~\ref{thm:exist} imply that 
$e(r) \rightarrow 0$ as $r \rightarrow 0$ and as $r \rightarrow \infty$.
Using the ODEs ((\ref{eq:ode1})-(\ref{eq:ode2})) we
can derive the equation
\[
  (-\Delta_r + \alpha)e + \frac{e}{f}e' = (1-\lambda)f^2f'
\]
where
\[
  \alpha(r) = \frac{1+n(1-a)}{r^2}(1 + \frac{rf'}{f}) 
  + f^2 + \frac{na'}{r} > 0
\]
and the result follows from the maximum principle.  $\Box$



\section{The linearized operator}
\label{sec:lin}

In this section, we introduce the linearized operator (or Hessian)
around the $n$-vortex, and identify its symmetry zero-modes.

\subsection{Definition of the linearized operator}
\label{section:linop}

We work on the real Hilbert space
\[
  X = L^2({\bf R}^2;{\bf C}) \oplus L^2({\bf R}^2;{\bf R}^2)
\]
with inner-product
\[
  < (\xi,B), (\eta,C) >_X = \int_{{\bf R}^2}
  \{ \Re (\bar{\xi}\eta) + B \cdot C \}.
\]
We define the linearized operator, $L_{\psi,A}$
(= the Hessian of $E(\psi,A)$) at a solution
$(\psi,A)$ of (\ref{eq:eq1})-(\ref{eq:eq2})
through the quadratic form
\[
  \frac{\partial^2}{\partial \epsilon \partial \delta}
  E(\psi + \epsilon\xi + \delta\eta, A + \epsilon B + \delta C)
  |_{\epsilon = \delta = 0} = 
  \langle (\eta, C) L_{\psi,A} (\xi, B) \rangle_X
\]
for all $(\xi, B)$, $(\eta, C)$, $\in X$.
The result is
\[
  L_{\psi,A} \left( \begin{array}{c} \xi \\ B \end{array} \right)
  = \left( \begin{array}{c}
  [-\Delta_A + \frac{\lambda}{2}(2|\psi|^2-1)]\xi + 
  \frac{\lambda}{2}\psi^2\bar{\xi}
  + i[2\nabla_A \psi + \psi\nabla]\cdot B  \\
  \Im( [\bar{\nabla_A\psi}-\bar{\psi}\nabla_A] \xi)
  + (-\Delta + \nabla\nabla + |\psi|^2) \cdot B  
  \end{array} \right).
\]


\subsection{Symmetry zero-modes}
\label{subsec:zeromodes}

We identify the part of the kernel of the operator
\[
  L^{(n)} \equiv L_{\psi^{(n)},A^{(n)}}
\]
which is due to the symmetry group.
\begin{prop}
\label{prop:modes}
We have
\begin{enumerate}
\item
\begin{equation}
\label{eq:gmode}
  L^{(n)} \left( \begin{array}{c} 
  i\gamma \psi^{(n)} \\ \nabla \gamma 
  \end{array} \right) = 0
\end{equation}
for any $\gamma : {\bf R}^2 \rightarrow {\bf R}$
\item
\begin{equation}
\label{eq:tmode}
  L^{(n)} \left( \begin{array}{c} 
  \partial_j \psi^{(n)} \\ \partial_j A^{(n)}
  \end{array} \right) = 0
\end{equation}
for $j=1,2$.
\end{enumerate}
\end{prop}
{\em Proof:}
We use the basic result that the generator of a one-parameter group
of symmetries of $E(\psi,A)$, applied to the $n$-vortex, lies in the kernel
of $L^{(n)}$.
The vector in~(\ref{eq:gmode}) is easily seen to be the generator
of a one-parameter family of gauge transformations~(\ref{eq:g1}-\ref{eq:g2}) 
applied to the $n$-vortex.  Similarly, the vector in~(\ref{eq:tmode})
is the generator of coordinate
translations applied to the $n$-vortex.  $\Box$
\begin{rem}
Applying the generator of the coordinate rotational symmetry~(\ref{eq:rot})
to the $n$-vortex gives us nothing new, it is contained in the 
gauge-symmetry case.
\end{rem}

We define $Z_{sym}$ to be the subspace of $X$ spanned by  
the $L^2$ zero-modes described in proposition~\ref{prop:modes}.
We recall that the $n$-vortex is called {\em stable} 
if there is a constant $c > 0$ such that
\begin{equation}
\label{eq:stabdef}
  L^{(n)} |_{Z_{sym}^{\perp}} \geq c,
\end{equation}
and {\em unstable} if $L^{(n)}$ has a negative eigenvalue.

\subsection{Gauge fixing}
\label{subsec:gauge}

In order to remove the infinite dimensional kernel of
$L^{(n)}$ arising from gauge symmetry, we restrict
the class of perturbations.  Specifically, we restrict
$L^{(n)}$ to the space of those perturbations
$(\xi, B) \in X$  
which are orthogonal to the $L^2$ gauge zero-modes~(\ref{eq:gmode}).
That is,
\[
      \left\langle \left( \begin{array}{c} 
      i\gamma \psi^{(n)} \\ \nabla \gamma 
      \end{array} \right), \left( \begin{array}{c} 
      \xi \\ B \end{array} \right) \right\rangle_X = 0
\] 
for all $\gamma$.  Integration by parts gives the gauge condition
\begin{equation}
\label{eq:gchoice}
  \Im(\overline{\psi^{(n)}}\xi) = \nabla \cdot B.
\end{equation}

As is done in \cite{s}, we consider a modified quadratic form
$\tilde{L}^{(n)}$, defined by
\[
  < \alpha, \tilde{L}^{(n)} \alpha > = 
  < \alpha, L^{(n)} \alpha >
  + \int (\Im(\overline{\psi^{(n)}} \xi) - \nabla \cdot B)^2
\]
for $\alpha = (\xi, B) \in X$.
Clearly, $\tilde{L}^{(n)}$ agrees with $L^{(n)}$ on the subspace
of $X$ specified by the gauge condition~(\ref{eq:gchoice}).
This modification has the important effect of shifting the essential
spectrum away from zero (see~(\ref{eq:contspec})).
A straightforward computation gives the following expression for 
$\tilde{L}^{(n)}$:
\[
  \tilde{L}^{(n)} \left( \begin{array}{c} \xi \\ B  \end{array} \right) = 
  \left( \begin{array}{c}
  [-\Delta_A + \frac{\lambda}{2}(2|\psi|^2-1) + \frac{1}{2}|\psi|^2]\xi
  + \frac{1}{2}(\lambda-1)\psi^2\bar{\xi} + 2i\nabla_A \psi \cdot B  \\
  2\Im[\bar{\nabla_A \psi}\xi] + [-\Delta + |\psi|^2]B 
  \end{array} \right).
\]
To establish theorem~\ref{thm:main},
it suffices to prove that $\tilde{L}^{(n)} \geq c > 0$
on the subspace of $X$
orthogonal to the translational zero-modes~(\ref{eq:tmode}).

$\tilde{L}^{(n)}$ is a real-linear operator on $X$.
It is convenient to identify $L^2({\bf R}^2;{\bf R}^2)$ with 
$L^2({\bf R}^2;{\bf C})$ through the correspondence
\begin{equation}
\label{eq:complex}
  B = \left( \begin{array}{c}  B_1 \\ B_2  \end{array} \right)
  \leftrightarrow B^c \equiv B_1 - iB_2,
\end{equation}
and then to complexify the 
space $X \mapsto \tilde{X} = [L^2({\bf R}^2;{\bf C})]^4$ via
\begin{equation}
\label{eq:comp}
  (\xi, B) \mapsto
  (\xi, \bar{\xi}, B^c, \bar{B}^c).
\end{equation}
As a result, $\tilde{L}^{(n)}$ is replaced by the 
complex-linear operator
\[
  \tilde{\tilde{L}}^{(n)} = 
  \mbox{ diag } \{ -\Delta_A, -\overline{\Delta_A}, -\Delta, -\Delta \}
  + V^{(n)}
\]
where
\[
  V^{(n)} = \left( \begin{array}{cccc}
  \frac{\lambda}{2}(2|\psi|^2-1) + \frac{1}{2}|\psi|^2  &
  \frac{1}{2}(\lambda-1)\psi^2  &  -i(\partial_A^* \psi)  &
  i(\partial_A \psi)  \\  \frac{1}{2}(\lambda - 1)\bar{\psi}^2  &
  \frac{\lambda}{2}(2|\psi|^2-1) + \frac{1}{2}|\psi|^2  & 
  -i(\bar{\partial_A \psi})  &  i(\bar{\partial_A^* \psi}) \\
  i(\bar{\partial_A^* \psi})  &  i(\partial_A \psi)  &
  |\psi|^2  &  0  \\   -i(\bar{\partial_A \psi})  &
  -i(\partial_A^* \psi)  &  0  &  |\psi|^2  \end{array} \right).
\]
Here we have used the notation
\[
  \partial_A \equiv \partial_z - iA
\]
where $\partial_z = \partial_1 - i\partial_2$ (and the superscript c has
been dropped from the complex function $A$ obtained from the vector-field
$A$ via~(\ref{eq:complex})).

The components of $V^{(n)}$ are bounded, and it follows from standard 
results (\cite{rs2}) that 
$\tilde{\tilde{L}}^{(n)}$ is a self-adjoint operator on $\tilde{X}$,
with domain
\[
  D(\tilde{\tilde{L}}^{(n)}) = [H_2({\bf R}^2;{\bf C})]^4
\]

\section{Block decomposition}
\label{sec:block}

We write functions on ${\bf R}^2$ in polar coordinates.  Precisely,
\begin{equation}
\label{eq:polar}
  \tilde{X} = [L^2({\bf R}^2; {\bf C})]^4 =
  [L^2_{rad} \otimes 
  L^2({\bf S}^1; {\bf C})]^4
\end{equation}
where $L^2_{rad} \equiv L^2({\bf R}^{+}, rdr)$.

Let $\rho_n : U(1) \rightarrow Aut([L^2({\bf S}^1; {\bf C})]^4)$ 
be the representation whose action is given by
\[
  \rho_n(e^{i\theta})
  (\xi, \eta, B, C)(x) = 
  (e^{in\theta}\xi, e^{-in\theta}\eta, e^{-i\theta}B,  
  e^{i\theta}C) (R_{-\theta}x)
\]
where $R_{\alpha}$ is a counter-clockwise rotation in ${\bf R}^2$
through the angle $\alpha$.
It is easily checked that the linearized operator $\tilde{\tilde{L}}^{(n)}$ 
commutes with  $\rho_n(g)$ for any $g \in U(1)$.  
It follows that $\tilde{\tilde{L}}^{(n)}$ leaves invariant
the eigenspaces of $d\rho_n(s)$ for any $s \in i{\bf R} = Lie(U(1))$.
The resulting block decomposition
of $\tilde{\tilde{L}}^{(n)}$, which is described in this section,
is essential to our analysis.  In particular, the translational zero-modes
each lie within a single subspace of this decomposition.

\subsection{The decomposition of $L^{(n)}$}
\label{section:blocks}

In what follows, we define, for convenience,
$b(r) = \frac{n(1-a(r))}{r}$.

\begin{prop}
\label{prop:decomp}
There is an orthogonal decomposition
\begin{equation}
\label{eq:orth}
  \tilde{X} = \bigoplus_{m \in {\bf Z}}
  ( e^{i(m+n)\theta} L^2_{rad} \oplus
  e^{i(m-n)\theta} L^2_{rad} \oplus
  -i e^{i(m-1)\theta} L^2_{rad} \oplus
  i e^{i(m+1)\theta} L^2_{rad} ),
\end{equation}
under which
the linearized operator around the vortex, $\tilde{\tilde{L}}^{(n)}$,
decomposes as
\[
  \tilde{\tilde{L}}^{(n)} = \bigoplus_{m \in {\bf Z}} \hat{L}_m^{(n)}
\]
where
\begin{equation}
  \label{eq:lm2}
  \hat{L}_m^{(n)} = -\Delta_r(Id) + \hat{V}_m^{(n)}
\end{equation}
with
\[
  \hat{V}_m^{(n)} = \frac{1}{r^2} 
  \mbox{ diag } \{ [m + n(1-a)]^2, [m - n(1-a)]^2, [m-1]^2, [m+1]^2 \}
  + V'
\]
and
\[
  V' = \left( \begin{array}{cccc} 
  \frac{\lambda}{2}(2f^2-1) + \frac{1}{2}f^2
  & \frac{1}{2}(\lambda - 1)f^2   &
  f' - bf  &  -[f' + bf]  \\ 
  \frac{1}{2}(\lambda-1)f^2  & 
  \frac{\lambda}{2}(2f^2-1) + \frac{1}{2}f^2  &  -[f' + bf]  
  &  f' - bf  \\  f' - bf  &
  -[f' + bf]  &  f^2  &  0  \\
  -[f' + bf]  &  f' - bf  &  0  &
  f^2  \end{array} \right).
\]
\end{prop}
{\em Proof:}
The decomposition~(\ref{eq:orth}) 
of $\tilde{X}$ follows from the usual Fourier decomposition
of $L^2({\bf S}^1;{\bf C})$, and the relation~(\ref{eq:polar}). 
An easy computation shows that
$\tilde{\tilde{L}}^{(n)}$
preserves the space of vectors of the form
\begin{equation}
\label{eq:eig}
  (\xi e^{i(m+n)\theta},
  \eta e^{i(m-n)\theta}, -i\alpha e^{i(m-1)\theta}, 
  i\beta e^{i(m+1)\theta} )
\end{equation}
and that it acts on such vectors via~(\ref{eq:lm2}).
$\Box$

It follows that $\hat{L}_m^{(n)}$ is self-adjoint on $[L^2_{rad}]^4$.
It will also be convenient to work with a rotated version
of the operator $\hat{L}_m^{(n)}$,
\[
  L_m^{(n)} \equiv \left\{ \begin{array}{cc}
            R \hat{L}_m^{(n)} R^T  &  m \geq 0  \\
            R' \hat{L}_m^{(n)} (R')^T  &  m < 0
            \end{array} \right. 
\]
where
\[
  R = \frac{1}{\sqrt{2}} \left( \begin{array}{cccc}
      1 & 1 & 0 & 0 \\
      -1 & 1 & 0 & 0 \\
      0 & 0 & 1 & 1 \\
      0 & 0 & 1 & -1  \end{array} \right),    \;\;\;\;\;\;\;\;\;\;
  R' = \frac{1}{\sqrt{2}} \left( \begin{array}{cccc}
       1 & 1 & 0 & 0 \\
       1 & -1 & 0 & 0 \\
       0 & 0 & 1 & 1 \\
       0 & 0 & 1 & -1  \end{array} \right).
\]
We have 
\begin{equation}
\label{eq:lm1}
  L_m^{(n)} = -\Delta_r(Id) + V_m^{(n)}
\end{equation}
where
\[
  V_m^{(n)} = \left( \begin{array}{cccc}
  \frac{m^2}{r^2} + b^2 + \frac{\lambda}{2}(3f^2-1)  &
  -2|m|\frac{b}{r}  &  -2bf  &  0  \\
  -2|m|\frac{b}{r}  &  
  \frac{m^2}{r^2} + b^2 + \frac{\lambda}{2}(f^2-1) + f^2  &
  0  &  -2f'  \\  -2bf  &  0  &  \frac{m^2+1}{r^2} + f^2  &
  -2\frac{|m|}{r^2}  \\  0  &  -2f'  & 
  -2\frac{|m|}{r^2}  &  \frac{m^2+1}{r^2} + f^2  \end{array} \right).
\]

\subsection{Properties of $L_m^{(n)}$}
\label{subsec:gen}

\begin{prop}
\label{prop:props}
We have the following:
\begin{enumerate}
\item
\begin{equation}
\label{eq:flip}
  L_m^{(n)} = L_{-m}^{(n)}
\end{equation}
\item
\begin{equation}
\label{eq:contspec}
  \sigma_{ess}(L_m^{(n)}) = [\min(1,\lambda), \infty)
\end{equation}
\item 
For $|n|=1$ and $m \geq 2$,
\begin{equation}
\label{eq:mon}
  L_m^{(n)} - L_1^{(n)} \geq 0
\end{equation}
with no zero-eigenvalue.
\end{enumerate}
\end{prop}
{\em Proof:}
The first statement is obvious.
The second statement follows in a standard way from the fact that 
\[
  \lim_{r \rightarrow \infty} V_m^{(n)}(r) =
  \mbox{ diag } \{\lambda, 1, 1, 1\}
\]
To prove the third statement, we compute
\[
  \hat{L}_m^{(n)} - \hat{L}_1^{(n)} = 
  \frac{m-1}{r^2} \mbox{ diag }
  \{ m+2n(1-a), \; m-2n(1-a), \; m-1, \; m+3 \}
\]
which is non-negative, with no zero-eigenvalue
for $m \geq 2$, $n=1$.  $\Box$

\begin{rem}
In light of~(\ref{eq:flip}), we can assume from now on that $m \geq 0$.
This degeneracy is a result of the complexification~(\ref{eq:comp})
of the space of perturbations.
\end{rem}

\subsection{Translational zero-modes}
\label{section:zeros}

The gauge fixing (section~\ref{subsec:gauge}) 
has eliminated the zero-modes arising from 
gauge symmetry.  The translational zero-modes remain.

As written in~(\ref{eq:tmode}), the translational zero-modes fail
to satisfy the gauge condition~(\ref{eq:gchoice}).  Further, they do
not lie in $L^2$.
A straightforward computation shows that 
if we adjust the vectors in~(\ref{eq:tmode}) by gauge zero-modes
given by~(\ref{eq:gmode}) with $\gamma = -A_j$, $j=1,2$, we obtain
\[
  T_1 = \left( \begin{array}{c}
  (\nabla_{A} \psi)_1  \\  
  (\nabla \times A) e_2
  \end{array} \right)
  \;\;\;\;\;\;
  T_2 = \left( \begin{array}{c}
  (\nabla_{A} \psi)_2  \\  
  -(\nabla \times A) e_1 
  \end{array} \right)
\]
where $e_1 = (1, 0)$ and $e_2 = (0, 1)$.
$T_1$ and $T_2$ satisfy~(\ref{eq:gchoice}), 
and are zero-modes of the 
linearized operator.  Note also that $T_{\pm 1}$ decay exponentially
as $|x| \rightarrow \infty$, and hence lie in $L^2$.

It is easily
checked that $T_1 \pm iT_2$ lie in the $m=\pm1$ blocks for $\hat{L}_m^{(n)}$.
After rotation by $R$, we have
 \[
  L_{\pm1}^{(n)} T = 0
\]
where
\[
  T = (f', b f, n\frac{a'}{r}, n\frac{a'}{r}).
\]

\section{Stability of the fundamental vortices}
\label{sec:fund}

In this section we prove the first part of theorem~\ref{thm:main}.
Specifically, we show that for some $c > 0$, $L_m^{(\pm 1)} \geq c$
for $m \not= 1$, and $L_1^{(\pm1)} |_{T^{\perp}} \geq c$.
In light of the discussions in sections~\ref{subsec:gauge}, 
\ref{section:blocks}, and \ref{section:zeros}, this will establish the 
stability of the $\pm1$-vortices.




\subsection{Non-negativity of $L_0^{(n)}$ and radial minimization}
\label{subsection:min}

\begin{prop}
\label{prop:rad}
  $L_0^{(n)} \geq 0$ for all $\lambda$.
\end{prop}
{\em Proof:}

From the expression~(\ref{eq:lm1}) we see that $L_0^{(n)}$
breaks up:
\begin{equation}
\label{eq:split}
  L_0^{(n)} = N_0 \oplus M_0
\end{equation}
(abusing notation slightly) where
\[
  M_0 = -\Delta_r(Id) + W_0
\]
with
\[
  W_0 = \left( \begin{array}{cc}
  b^2 + \frac{\lambda}{2}(3f_n^2-1) & -2bf \\
  -2bf & \frac{1}{r^2} + f^2
  \end{array} \right)
\]
and
\[
  N_0 = \left( \begin{array}{cc} 
  -\Delta_r + b^2 + \frac{\lambda}{2}(f^2-1) + f^2  &
  -2f'  \\  -2f'  &  -\Delta_r + \frac{1}{r^2} + f^2.
  \end{array} \right)
\]
An easy computation shows that
$M_0$ is precisely the Hessian of the radial energy, 
$Hess E_r^{(n)}$~(see (\ref{eq:erad})).
Since the $n$-vortex minimizes $E_r^{(n)}$,
we have $M_0 \geq 0$.  It remains to show $N_0 \geq 0$.  We establish
the stronger result, $N_0 > 0$.  Note that
\[
  N_0 = G_0^* G_0
\]
where
\[
  G_0 = \left( \begin{array}{cc}
  \partial_r - f'/f  &  f  \\
  f  &  \partial_r + 1/r 
  \end{array} \right)
\]
In fact, $G_0$ has no zero-eigenvalue.  To see this, we first
remark that $G_0$ is a relatively compact perturbation of
$G_0|_{\lambda = 1}$, due to the exponential decay of the field
components.  It follows from an index-theoretic calculation done 
in~\cite{w,s}, that $G_0|_{\lambda = 1}$ is Fredholm, with index $0$.
We conclude that the same is true of $G_0$ (for any $\lambda$).  
Finally, it is a simple matter to check that $G_0^*$ has trivial kernel.  If
\[
  G_0^* \left( \begin{array}{c} \xi \\ \beta \end{array} \right) = 0
\]
it follows that
\[
  (-\Delta_r + f^2)\beta = 0
\]
and hence that $\beta = 0$, and so $\xi = 0$.
The relation $N_0 > 0$ follows from this, and the the fact that
$\sigma_{ess}(N_0) = [1,\infty)$.  $\Box$

\subsection{A maximum principle argument}

Removing the equality in proposition~\ref{prop:rad} requires more work. 
First, we establish an extension of the maximum
principle to systems (see, eg, \cite{lm,pa} for related results).  
We will use this
also in the proof that the the translational zero-mode is the
ground state of $L_1^{(n)}$ (section~\ref{subsec:L1}).

\begin{prop}
\label{prop:mainmp}
Let $L$ be a self-adjoint operator on $L^2({\bf R}^n;{\bf R}^d)$ of the form
\[
  L = -\Delta (Id) + V
\]
where $V$ is a $d \times d$ matrix-multiplication operator with 
smooth entries.  Suppose that $L \geq 0$ and that 
for $i \not= j$, $V_{ij}(x) \leq 0$ for all $x$.  
Further, suppose $V$ is irreducible in the sense that for any splitting
of the set $\{1,\ldots,d\}$ into disjoint sets $S_1$ and $S_2$, there is
an $i \in S_1$ and a $j \in S_2$ with $V_{ij}(x) < 0$ for all $x$.  
Finally, suppose that
$L \xi = \eta \in L^2$ with $\eta \geq 0$ component-wise, and
$\xi \not\equiv 0$.
Then either
\begin{enumerate}
\item $\xi > 0$ or
\item $\eta \equiv 0$ and $\xi < 0$.
\end{enumerate}
\end{prop} 
{\em Proof:}
We write $\xi = \xi^+ - \xi^-$ with $\xi^+, \xi^- \geq 0$ component-wise,
and compute 
\[
  0 \;\; \leq \;\; <\xi^-, L \xi^-> \;\; = \;\; <\xi^-, L \xi^+> 
  - <\xi^-, L \xi>.
\]
Since $\xi_j^+$ and $\xi_j^-$ have disjoint support, we have
\[
  r.h.s = \sum_{j \not= k} <\xi_j^-, V_{jk} \xi_k^+> 
  - <\xi^-,\eta> \;\; \leq 0.
\]
Thus we have
\begin{enumerate}
\item $0 = \;\; <\xi^-, L\xi^->$
\item $0 = \;\; <\xi_j^-, V_{jk} \xi_k^+>$ for all $j \not= k$
\end{enumerate}
Since $L \geq 0$, the first of these implies $L \xi^- = 0$
and hence $L \xi^+ = \eta$.  So if $\eta \not\equiv 0$, then
$\xi^+ \not\equiv 0$.  If $\eta \equiv 0$ and $\xi^+ \equiv 0$,
replace $\xi$ with $-\xi$ in what follows.
An application of the strong maximum principle 
(eg. \cite{gt}, Thm. 8.19) to each component of the equation
\[
  L \xi^+ = \eta
\]
now allows us to conclude that for each $k$, either $\xi_k^+ > 0$
or $\xi_k^+ \equiv 0$.
We know that for some $k$, $\xi_k^+ > 0$.
Looking back at the second listed equation
above, and using the irreducibility of $V$, 
we then see that $\xi_j^- \equiv 0$ for all $j$.
Finally, we can easily rule out the possibility $\xi_k \equiv 0$ 
for some $k$, by looking back at the equation satisfied by $\xi_k$.
Thus we have $\xi > 0$.  $\Box$

\subsection{Positivity of $L_0^{(n)}$}
\label{subsection:l0}

Now we apply proposition~\ref{prop:mainmp} to show
$M_0 > 0$.  The trick here is to find a function $\xi$ which satisfies
$M_0 \xi \geq 0$.  This allows us to rule out the existence of a 
zero-eigenvector, which would be positive by proposition~\ref{prop:mainmp}.
To obtain such a $\xi$, we differentiate the vortex with respect to the
parameter $\lambda$.  Specifically,  
differentiation of the Ginzburg-Landau 
equations with respect to $\lambda$ results in
\begin{equation}
\label{eq:diff}
  M_0 \xi = \eta
\end{equation}
where
\[
  \xi = \left( \begin{array}{cc}
  \partial_{\lambda} f  \\  n \partial_{\lambda} a/r
  \end{array} \right)
\]
and
\[
  \eta = \left( \begin{array}{cc} 
  \frac{1}{2}(1-f^2)f  \\  0  \end{array} \right) \geq 0.
\]
We can now establish
\begin{prop}
\label{prop:l02}
For all $\lambda$, $L_0^{(n)} \geq c > 0$.
\end{prop}
{\em Proof:}
We have already shown in the proof of proposition~\ref{prop:rad}, that
$N_0 > 0$ and $M_0 \geq 0$.  Hence, due to~(\ref{eq:split})
and~(\ref{eq:contspec}),
it suffices to show that $Null(M_0) = \{0\}$.
Suppose $M_0 \zeta = 0$, $\zeta \not\equiv 0$.
Proposition~\ref{prop:mainmp} then implies $\zeta > 0$
(or else take $-\zeta$).  Now
\[
  0 = \;\; <M_0\zeta, \xi> \;\; = \;\; <\zeta, M_0\xi> \;\; 
  = \;\; <\zeta, \eta> \;\; > \;\; 0
\]
gives a contradiction.  $\Box$
\begin{rem}
Proposition~\ref{prop:mainmp} applied to equation~(\ref{eq:diff}) also
gives $\xi > 0$.  That is, the vortex profiles increase 
monotonically with $\lambda$.  This can be used to show that the rescaled
vortex $(f_n(r/\sqrt{\lambda}), a_n(r/\sqrt{\lambda}))$ converges as
$\lambda \rightarrow \infty$ to $(f^*, 0)$, where $f^*$ is the (profile
of) the $n$-vortex solution of the ordinary GL equation:
$-\Delta_r f^* + n^2 f^*/r^2 + ({f^*}^2-1)f^* = 0$.  This result was 
established by different means in~\cite{abg}.
\end{rem}

\subsection{Positivity of $L_1^{(\pm1)}$}
\label{subsec:L1}

\begin{prop}
\label{prop:l1}
$ L_1^{(\pm1)} \geq 0 $
with non-degenerate zero-eigenvalue given by $T$.
\end{prop}

{\em Proof:}
Let $\mu = inf spec L_1^{(\pm 1)} \leq 0$, which is an 
eigenvalue by~(\ref{eq:contspec}).
Suppose $L_1^{(\pm1)} S = \mu S$. Applying proposition~\ref{prop:mainmp} 
to $L_1^{(\pm 1)} - \mu$ (note that $V_1^1$ satisfies the irreducibility
requirement) gives  
$S > 0$ (or $S < 0$).  Further, $\mu$ is
non-degenerate, as if $\mu$ were degenerate, we would have two 
strictly positive eigenfunctions which are orthogonal, an impossibility.
Now if $\mu < 0$, we have
$<S,T> = 0$, which is also impossible.  Thus $S$ is a multiple of $T$, 
and $\mu = 0$.  $\Box$

\subsection{Completion of stability proof for $n=\pm1$}
\label{subsec:mon}

We are now in a position to complete the proof of the first 
statement of theorem~\ref{thm:main}.
By proposition~\ref{prop:l02},
$L_0^{(\pm1)} \geq c > 0$.
By proposition~\ref{prop:l1} and~(\ref{eq:contspec}), 
$L_1^{(\pm1)}|_{T^{\perp}} \geq \tilde{c} > 0$.
Finally, by~(\ref{eq:mon}), $L_m^{(\pm1)} \geq c' > 0$
for $|m| \geq 2$.  It follows from proposition~\ref{prop:decomp} that
$\tilde{L}^{(n)} \geq c > 0$ on the subspace of $X$ orthogonal
to the translational zero-modes.  By the discussion of 
section~\ref{subsec:gauge}, this gives theorem~\ref{thm:main} for $n = \pm1$.
$\Box$

\section{The critical case, $\lambda = 1$}
\label{sec:crit}

In order to prove the remainder of theorem~\ref{thm:main},
we exploit some results from the $\lambda = 1$ case.

\subsection{The first-order equations}

Following \cite{bog}, we use an integration by parts
to rewrite the energy~(\ref{eq:ac}) as
\begin{equation}
  \label{eq:bog}
  E(\psi,A) = \frac{1}{2} \int_{{\bf R}^2} \{ |\partial_A \psi|^2
  + [\nabla \times A + \frac{1}{2}(|\psi|^2-1)]^2 +
  \frac{1}{4}(\lambda - 1)(|\psi|^2-1)^2 \}
  + \pi \mbox{deg}(\psi)
\end{equation}
(recall, since we work in dimension two, $\nabla \times A$ is a scalar)
where $\mbox{deg}(\psi)$ is the topological degree of $\psi$,
defined in the introduction.
We assume, without loss of generality, that $deg(\psi) \geq 0$.
Clearly, when $\lambda = 1$, a solution of the first-order equations
\begin{equation}
  \label{eq:first1}
  \partial_A \psi = 0
\end{equation}
\begin{equation}
  \label{eq:first2}
  \nabla \times A + \frac{1}{2}(|\psi|^2-1) = 0
\end{equation}
minimizes the energy within a fixed topological
sector, $\deg(\psi) = n$, and hence is stable.
Note that we have identified the vector-field $A$ with a complex
field as in~(\ref{eq:complex}).

The $n$-vortices~(\ref{eq:nvortex}) 
are solutions of these equations (when $\lambda = 1$).
Specifically,
\begin{equation}
\label{eq:firstode1}
  n\frac{a'}{r} = \frac{1}{2}(1-f^2)
\end{equation}
and
\begin{equation}
\label{eq:firstode2}
  f' = n\frac{(1-a)f}{r}.
\end{equation}
In fact, it is shown in~\cite{t2} that for $\lambda = 1$,
any solution of the variational equations solves the first-
order equations~(\ref{eq:first1}-\ref{eq:first2}).

Beginning from expression~(\ref{eq:bog}) 
for the energy, the variational
equations (previously written 
as~(\ref{eq:eq1}-\ref{eq:eq2}))
can be written as
\begin{equation}
  \label{eq:eq21}
  \partial_A^*[\partial_A \psi] + \psi[\nabla \times A
  + \frac{1}{2}(|\psi|^2-1)] + \frac{1}{2}(\lambda - 1)
  (|\psi|^2-1)\psi = 0
\end{equation}
\begin{equation}
  \label{eq:eq22}
  i\bar{\psi}[\partial_A \psi] - i\partial_{\bar{z}}[\nabla \times A
  + \frac{1}{2}(|\psi|^2-1)] = 0
\end{equation}
(here $\partial_A^* \equiv -\partial_z + iA$ is the
adjoint of $\partial_A$).

\subsection{First-order linearized operator}

We show that the linearized operator at $\lambda = 1$ is the 
square of the linearized operator for the first-order equations.

Linearizing the first-order equations~(\ref{eq:first1}-\ref{eq:first2})
about a solution,
$(\psi, A)$ (of the first-order equations) 
results in the following equations for the perturbation,
$\alpha \equiv (\xi, B)$:
\[
  \partial_A\xi - iB\psi = 0
\]
\[
  \nabla \times B + \Re(\bar{\psi}\xi) = 0.
\]
Now using $-i\partial_zB = \nabla \times B - i(\nabla \cdot B)$,
and adding in the gauge condition~(\ref{eq:gchoice}),
we can rewrite this as
\begin{equation}
  \label{eq:1lin}
  L_1 \alpha = 0
\end{equation}
where 
\[
  L_1 = \left( \begin{array}{cc}
  \partial_A  &  -i\psi  \\
  \bar{\psi}   &  -i\partial_z 
  \end{array} \right).
\]


If we linearize the full (second order) variational
equations (in the form~(\ref{eq:eq21}-\ref{eq:eq22}))
around $(\psi, A)$, we obtain
\[
  \partial_A^*[\partial_A\xi - iB\psi] + i\bar{B}[\partial_A \psi]
  + \psi[\nabla \times B + \Re(\bar{\psi}\xi)] 
\]
\[
  + \xi[\nabla \times A + \frac{1}{2}(|\psi|^2-1)]
  + \frac{1}{2}(\lambda - 1)[(|\psi|^2-1)\xi + 2\psi\Re(\bar{\psi}\xi)] = 0
\]
and
\[
  i\bar{\psi}[\partial_A\xi - iB\psi] + i\bar{\xi}[\partial_A \psi]
  -i\partial_{\bar{z}}[\nabla \times B + \Re(\bar{\psi}\xi)] = 0.
\]
\begin{prop}
When $\lambda = 1$,
these linearized equations can also be written
\[
  L_1^* L_1 \alpha = 0
\]
\end{prop}
{\em Proof:}
This is a simple computation using the fact that
the first-order equations~(\ref{eq:first1}-\ref{eq:first2}) hold.  $\Box$

This relation holds also on the level of the blocks.
A straightforward computation gives
\[
  L_m^{(n)}|_{\lambda = 1} = F_m^* F_m
\]
where
\[
  F_m = \left( \begin{array}{cccc}
        \partial_r - b  &  \frac{m}{r}  &  0  &  f  \\
        \frac{m}{r}  &  \partial_r - b  &  -f  &  0  \\
        0  &  -f  &  \partial_r + 1/r  &  \frac{m}{r}  \\
        f  &  0  &  \frac{m}{r}  &  \partial_r + 1/r
        \end{array} \right)
\]

\subsection{Zero-modes for $\lambda = 1$}

It was predicted in \cite{w} (and proved rigorously in \cite{s}) that
for $\lambda = 1$, the linearized operator around any degree-$n$
solution of the first-order equations has a $2|n|$-dimensional
kernel (modulo gauge transformations).  This kernel arises because the
Taubes solutions form a $2|n|$-parameter family, and all have the
same energy.  The zero-eigenvalues
are identified in \cite{bog}, and we describe them here.
Let $\chi_m$ be the unique solution of
\[
  (-\Delta_r + \frac{m^2}{r^2} + f^2) \chi_m = 0
\]
on $(0,\infty)$ with
\[
  \chi_m \sim r^{-m} \;\;\;\; \mbox{ as } \;\;\;\;  r \rightarrow 0
\]
and
\[
  \chi_m \rightarrow 0  \;\;\;\; \mbox{ as } \;\;\;\; r \rightarrow \infty
\]
for $m = 1,2,\ldots,n$.  Then it is easy to check that
\begin{equation}
\label{eq:bogmodes}
  F_{\pm m} W_m = 0
\end{equation}
where 
\[
  W_m = \left( \begin{array}{c}
  f\chi_m  \\  f\chi_m  \\  -(\chi_m' + m\chi_m/r)  \\  
  -(\chi_m' + m\chi_m/r)
  \end{array} \right).
\]
We remark that
\[
  \chi_1 = \frac{1-a}{r}
\]
and it is easily verified that for $\lambda = 1$,
$W_{\pm 1} = T$ are the translational zero-modes.

\section{The (in)stability proof for $|n| \geq 2$}
\label{sec:high}

Here we complete the proof of theorem~\ref{thm:main}.

The idea is to decompose $L_m^{(n)}$ into a sum of two
terms, each of which has the same (translational) zero-mode 
(for $m=1$) as $L_m^{(n)}$.  One term is manifestly
positive, and the other satisfies restrictions of 
Perron-Frobenius theory.

We begin by modifying $F_m$, and defining, for any $\lambda$,
\[
  \tilde{F}_m  \equiv 
  \left( \begin{array}{cccc}
  (\partial_r - \frac{f'}{f}) \cdot q  &  \frac{m}{r}  &  0  &  f  \\
  \frac{m}{r}q  &  \partial_r  - \frac{f'}{f}  &  -f  &  0  \\
  0  &  -f  &  \partial_r + 1/r  &  \frac{m}{r}  \\
  fq  &  0  &  -\frac{m}{r}  &  \partial_r + 1/r
  \end{array} \right)
\]
where we have defined
\begin{equation}
\label{eq:q}
  q(r) \equiv \frac{n(1-a)f}{rf'}
\end{equation}
and $\partial_r \cdot q$ denotes an operator composition.
By~(\ref{eq:firstode2}), we have
$q \equiv 1$ for $\lambda = 1$.
We also set, for $m=1, \ldots, n$,
\[
  \tilde{W}_m = \left( \begin{array}{c}
  q^{-1} f \chi_m  \\  f \chi_m  \\  
  -( \chi_m' + m \frac{\chi_m}{r} ) \\
  -( \chi_m' + m \frac{\chi_m}{r} ) 
  \end{array} \right)
\]
Now $\tilde{W_m}$ has the following properties:
\begin{enumerate}
\item $\tilde{W}_{\pm 1}$ is the translational zero-mode
      $T$ for all $\lambda$
\item when $\lambda = 1$, $\tilde{W}_m = W_m$, $m = \pm 1, \ldots, \pm n$,
      give the $2n$ zero-modes~(\ref{eq:bogmodes}) of the linearized operator.
\end{enumerate}
These $W_m$ were chosen in \cite{bog} as candidates 
for directions of energy decrease (for $|m| \geq 2$) 
when $\lambda > 1$.
Intuitively, we think of $\tilde{W}_m$ as a perturbation that tends
to break the $n$-vortex into separate vortices of lower degree.

Now, $\tilde{F}_m$ was designed to have the following properties:
\begin{enumerate}
\item $\tilde{F}_m = F_m$ when $\lambda = 1$ (this is clear)
\item $\tilde{F}_m \tilde{W}_m = 0$ for all $m$ and $\lambda$
      (this is easily checked).
\end{enumerate}

A straightforward computation gives
\begin{equation}
\label{eq:keysplit}
  L_m^{(n)} = \tilde{F}_m^* \tilde{F}_m
  + J M_m
\end{equation}
where $J = \mbox{diag} \{ 1, 0, 0, 0 \}$ and
\[
  M_m = l_m - ql_mq + (\lambda - q^2)f^2
\]
with
\[
  l_m = -\Delta_r + \frac{m^2}{r^2} + b^2 + \frac{\lambda}{2}(f^2-1).
\]
By construction, when $m=1$, the second term in the
decomposition~(\ref{eq:keysplit}) must have a zero-mode corresponding to the
original translational zero-mode.  In fact, one can easily
check that $M_1 f' = 0$.

\begin{prop}
\label{prop:mp}
For $|n| \geq 2$,
$M_1$ has a non-degenerate zero-eigenvalue corresponding to $f'$,
and
\[
  \left \{  \begin{array}{cc}
  M_1 \geq 0  &  \lambda < 1  \\
  M_1 \leq 0  &  \lambda > 1
  \end{array} \right.
\]
on $L^2_{rad}$.
\end{prop}
{\em Proof:}
We recall inequality~(\ref{eq:ineq}), which implies that
for $\lambda < 1$, $q < 1$, and for $\lambda > 1$, $q > 1$. 
The operator $M_1$ is of the form
\begin{equation}
\label{eq:firstterm}
  M_1 = (1-q^2)(-\Delta_r) + \mbox{ first order } + \mbox{ multiplication. }
\end{equation}
One can show that
$M_1$ is bounded from below (resp. above) for $\lambda < 1$ 
(resp. $\lambda > 1$).
We stick with the case $\lambda < 1$ for concreteness.
Suppose $M_1 \eta = \mu \eta$ with $\mu = infspec M_1 \leq 0$.
Applying the maximum principle (eg proposition~\ref{prop:mainmp} for $d=1$)
to~(\ref{eq:firstterm}), we conclude 
that $\eta > 0$.  If $\mu < 0$, we have
$<\eta, f'> = 0$, a contradiction.  Thus $\mu = 0$, and is non-degenerate
by a similar argument.
$\Box$

We also have
\begin{lemma}
\label{lemma:mon}
For $m \geq 2$, $M_m - M_1$
is non-negative for $\lambda < 1$, non-positive for $\lambda > 1$,
and has no zero-eigenvalue.
\end{lemma}
{\em Proof:}  This follows from the equation
\[
  M_m - M_1 = (1-q^2)\frac{m^2-1}{r^2}.  \;\;\;\;\;\;\;\;  \Box
\]

{\em Completion of the proof of theorem~\ref{thm:main}:}
Suppose now $\lambda < 1$. 
Since $\tilde{F}_m^* \tilde{F}_m$ is manifestly non-negative,
and $M_m > M_1$ for $m \geq 2$, we have
$L_m^{(n)} \geq 0$ for $m \geq 1$ (with only the
translational $0$-mode).  Combined with~(\ref{eq:contspec}) and 
propositions~\ref{prop:l02} and~\ref{prop:decomp}, this gives stability
of the $n$-vortex for $\lambda < 1$.

Now suppose $\lambda > 1$.  By~(\ref{eq:keysplit}), proposition~\ref{prop:mp}
and lemma~\ref{lemma:mon},
we have for $m = \pm 2, \ldots \pm n$,
\[
  <\tilde{W}_m, L_m^{(n)} \tilde{W}_m > \;\;\; < \;\; 0.
\]
We remark that $\tilde{W}_m$ corresponds to an element of the un-complexified
space $X$, and so $L^{(n)}$ has negative eigenvalues.
This establishes the instability of the $n$-vortex 
for $|n| \geq 2$, $\lambda > 1$,
and completes the proof of theorem~\ref{thm:main}.
$\Box$

\section{Appendix: vortex solutions are radial minimizers}

\begin{prop}
For $\lambda \geq 2n^2$,
a solution of the equations~(\ref{eq:ode1}-\ref{eq:ode2}) 
minimizes $E_r^{(n)}$.
\end{prop}
{\em Proof:}
It suffices then to show $M_0 = Hess E_r^{(n)} > 0$ (see 
section~\ref{subsection:min}).
We write $M_0 = L_0 + Z_0$ where
\[
  L_0 = diag\{ l, -\Delta_r \}
\]
with $l = -\Delta_r + b^2 + \frac{\lambda}{2}(f^2-1)$ and
\[
  Z_0 = \left( \begin{array}{cc}
        2\lambda f^2  &  -2bf  \\
	-2bf  &  \frac{1}{r^2} + f^2   \end{array} \right).
\]
We note that $l f = 0$ (one of the GL equations).  It follows from the 
fact that $f > 0$ and a Perron-Frobenius type argument (see \cite{os1})
that $l \geq 0$ with no zero-eigenvalue.  It suffices to show $Z_0 \geq 0$.
Clearly $tr(Z_0) > 0$, and
\[
  \det(Z_0) = 2\lambda f^4 + \frac{2f^2}{r^2}[\lambda - 2n^2(1-a)^2]
\]
is strictly positive for $\lambda \geq 2n^2$.  $\Box$


\begin{thebibliography}{9999}

\bibitem[ABG]{abg} L. Almeida, F. Bethuel, Y. Guo: A remark on the 
instability of symmetric vortices with large coupling constant.
Commun. Pure Appl. Math. {\bf 50} (1997) 1295-1300.

\bibitem[BC]{bc} M. S. Berger, Y. Y. Chen:  Symmetric vortices for the
nonlinear Ginzburg-Landau equations of superconductivity, and the
nonlinear desingularization phenomenon. J. Fun. Anal. {\bf 82}
(1989) 259-295.

\bibitem[B]{bog} E. B. Bogomol'nyi:  The stability of classical solutions.
Yad. Fiz. {\bf 24} (1976) 861-870.

\bibitem[BGP]{bgp} A. Boutet de Monvel-Berthier, V. Georgescu, R. Purice:
A boundary value problem related to the Ginzburg-Landau model.
Comm. Math. Phys. {\bf 142} (1991) 1-23.

\bibitem[GT]{gt} D. Gilbarg, N.S. Trudinger: Elliptic Partial 
Differential Equations of Second Order.  Berlin: Springer-Verlag, 1977.

\bibitem[G1]{g1} S. Gustafson: Symmetric solutions of Ginzburg-Landau
equations in all dimensions.  Intern. Math. Res. Notices (1997) No. 16, 
807-816.

\bibitem[G2]{g2} S. Gustafson:  Dynamical stability of magnetic vortices.
In preparation.

\bibitem[JT]{jt} A. Jaffe, C. Taubes:  Vortices and Monopoles.
Boston:  Birkhauser, 1980.

\bibitem[JR]{jr} L. Jacobs, C. Rebbi:  Interaction of superconducting
vortices.  Phys. Rev. {\bf B19} (1979) 4486-4494.

\bibitem[LM]{lm} J. Lopez-Gomez, M. Molina-Meyer:  The maximum principle
for cooperative weakly coupled elliptic systems and some applications.
Diff. Int. Eqns., 7 (1994) no. 2, 383-398.

\bibitem[OS1]{os1} Y. Ovchinnikov, I.M. Sigal:  Ginzburg-Landau
equation I:  Static Vortices. PDEs and their Applications,
Greiner et. al., eds. Providence:  AMS, 1997.

\bibitem[P]{p} B. Plohr:  Princeton thesis.

\bibitem[PA]{pa} C.V. Pao:  Nonlinear elliptic systems in unbounded
domains.  Nonlinear Analysis:  Theory, Methods, and Applications
{\bf 22} (1994) No. 11, 1391-1407.

\bibitem[RSII]{rs2} M. Reed, B. Simon:  Methods of Modern Mathematical
Physics, Vol II.  New York:  Academic Press, 1972.

\bibitem[RSIV]{rs4} Reed and Simon, Vol IV. 

\bibitem[S]{s} D. Stuart:  Dynamics of Abelian Higgs vortices in 
the near Bogomolny regime.  Commun. Math. Phys. {\bf 159} (1994) 51-91.

\bibitem[T1]{t1} C. Taubes:  Arbitrary $n$-vortex solutions to the 
first order Ginzburg-Landau equations.  Commun. Math. Phys. 
{\bf 72} (1980) 277.

\bibitem[T2]{t2} C. Taubes:  On the equivalence of the first and second
order equations for gauge theories.  Commun. Math. Phys. {\bf 75} (1980) 207.

\bibitem[W]{w} E. Weinberg:  Multivortex solutions of the Ginzburg-Landau 
equations.  Phys. Rev. D {\bf 19} (1979) 3008-3012.

\end{thebibliography}
\end{document}